\documentclass[11pt,twoside]{article}
\usepackage{amsmath}
\usepackage{amsthm}
\usepackage{amssymb}
\usepackage{amsfonts}

\input xypic

\makeatletter

\newtheorem{The}{Theorem}[section]

\newtheorem{Pro}[The]{Proposition}

\newtheorem{Def}[The]{Definition}

\newtheorem{Conj}[The]{Conjecture}

\def\endproof{\relax\ifmmode\expandafter\endproofmath\else
\unskip\nobreak\hfil\penalty50\hskip.75em\hbox{}\nobreak\hfil\bull
{\parfillskip=0pt \finalhyphendemerits=0 \bigbreak}\fi}
\def\endproofmath$${\eqno\bull$$\bigbreak}
\def\bull{\vbox{\hrule\hbox{\vrule\kern3pt\vbox{\kern6pt}\kern3pt\vrule}\hrule}}
 
\def\cancel#1#2{\ooalign{$\hfil#1\mkern1mu/\hfil$\crcr$#1#2$}}

\def\dirac{\mathpalette\cancel\partial}
 
\def\gr{{\rm gr}}
\def\cF{{\mathcal F}}

\newcommand{\ba}{\begin{eqnarray}}
\newcommand{\na}{\end{eqnarray}}
\newcommand{\ban}{\begin{eqnarray*}}
\newcommand{\nan}{\end{eqnarray*}}
\newcommand{\spinc}{\mathrm{Spin}^c}

\newcommand{\scr}{\mathcal}
\newcommand{\C}{\mathbb{C}}

\newcommand{\R}{\mathbb{R}}
\newcommand{\Z}{\mathbb{Z}}

\newcommand{\M}{\scr{M}}
\newcommand{\F}{\scr{F}}
\newcommand{\la}{\langle}
\newcommand{\ra}{\rangle}

\newcommand{\s}{\mathfrak{s}}

\newcommand{\disp}{\displaystyle}

\begin{document}
 
\title{Variants of Equivariant Seiberg-Witten Floer Homology\\
}
\author{Matilde Marcolli,  Bai-Ling Wang }
\date{}
\maketitle
 
\begin{abstract}
For a rational homology 3-sphere $Y$ with a $\spinc$ structure $\s$,
we show that simple algebraic manipulations of our construction of
equivariant Seiberg-Witten Floer homology in \cite{MW} lead to a 
collection of variants $HF^{SW, -}_{*,U(1)}(Y, \s)$, 
$HF^{SW, \infty}_{*, U(1)}(Y, \s)$
$HF^{SW, +}_{*, U(1)}(Y, \s)$, $\widehat{HF}^{SW}_{*}(Y, \s)$
and $HF^{SW}_{red, *}(Y, \s)$ which are topological invariants.
We establish a long exact sequence relating  
$HF^{SW, \pm}_{*, U(1)}(Y, \s)$ and $HF^{SW, \infty}_{*, U(1)}(Y, \s)$. 
We show they satisfy a duality under orientation reversal, and we 
explain their relation to the equivariant Seiberg-Witten Floer 
(co)homologies introduced in \cite{MW}. We conjecture
the equivalence of these versions of equivariant Seiberg-Witten
Floer homology with the Heegaard Floer invariants introduced by
Ozsv\'ath and Szab\'o.

\noindent{\bf Key words:} rational homology 3-spheres, equivariant
Seiberg-Witten Floer homology, $\spinc$ structures, topological
invariants. \par
\noindent{\bf Mathematics Subject Classification.} Primary:
57R58. Secondary: 57R57, 58J10. 

\end{abstract}

\section{Introduction}

For any rational homology 3-sphere $Y$ with a $\spinc$ structure $\s$, we 
constructed in \cite{MW} an equivariant Seiberg-Witten Floer homology 
$HF^{SW}_{*, U(1)}(Y, \s)$, which is a topological invariant.  In this paper,
we will generalize this construction to provide a collection of equivariant
 Seiberg-Witten Floer homologies $HF^{SW, -}_{*, U(1)}(Y, \s),
HF^{SW, \infty}_{*, U(1)}(Y, \s), 
HF^{SW, +}_{*, U(1)}(Y, \s), \widehat{HF}^{SW}_{*}(Y, \s)$
and $HF^{SW}_{red, *}(Y, \s)$, all of which are topological invariants,
such that $HF^{SW, +}_{*, U(1)}(Y, \s)$ is isomorphic to the
 equivariant Seiberg-Witten Floer homology $HF^{SW}_{*, U(1)}(Y, \s)$
constructed in \cite{MW}. 
The construction utilizes the $U(1)$-invariant
forms on $U(1)$-manifolds twisted with coefficients in the Laurent
polynomial algebra 
over integers.

In analogy to Austin and Braam's construction of equivariant instanton Floer
homology in \cite{AB}, the equivariant Seiberg-Witten Floer homology
$HF^{SW}_{*, U(1)}(Y, \s)$ is the homology of the complex
$(CF^{SW}_{*, U(1)}(Y, \s), D)$, where $CF^{SW}_{*, U(1)}(Y, \s)$
is generated by equivariant de Rham forms over all 
$U(1)$-orbits of the solutions of 3-dimensional Seiberg-Witten 
equations on $(Y, \s)$ modulo based gauge transformations (Cf.\cite{MW}). 
More specifically, 
\ba
CF^{SW}_{*, U(1)}(Y, \s) = \bigoplus_{a\in \M^{*}_Y(\s)}
\Z[\Omega] \otimes (\Z\eta_a \oplus \Z 1_a) \oplus 
\Z[\Omega] \otimes \Z 1_\theta,
\label{CF+}
\na
where $\M_Y(\s) = \M^{*}_Y(\s) \cup \{ \theta\}$ is the equivalence classes of
solutions to the Seiberg-Witten equations for a good pair of metric and
perturbations, consists of the irreducible monopoles $\M^{*}_Y(\s)$ and
the unique reducible monopole $\theta$. We used the notation $\eta_a$
to denote a 1-form on $O_a\cong S^1$, such that the cohomology class
$[\eta_a]$ is an integral generator of $H^1(O_a)$. Similarly, we
denote by $1_a$ the 0-form given by the constant function.

Each generator is endowed with
a grading such that, for any $k\geq 0$,
\ba\label{grade}
gr(\Omega^k\otimes \eta_a ) = 2k + gr (a), \ 
gr(\Omega^k\otimes 1_a) = 2k +gr (a) +1, \
and \  \ gr (\Omega^k\otimes 1_\theta) = 2k,
\na
where $gr:  \M^{*}_Y(\s) \to \Z$ is the relative grading with respect
to the reducible monopole $\theta$. This corresponds to
grading equivariant de Rham forms on each orbit $O_a$ by codimension
(Cf.\cite{MW} \S 5 for details).

The differential
operator $D$ can be expressed explicitly in components as the form:
\ba
\begin{array}{lll}
D(\Omega^k\otimes \eta_a ) &=& 
\disp{\sum_{b\in \M^{*}_Y(\s)\atop gr(a)-gr(b)=1} }n_{ab}
\Omega^k\otimes \eta_b + \sum_{c\in \M^{*}_Y(\s)\atop gr(a)-gr(c)=2}
m_{ac}\Omega^k\otimes 1_c
 -\Omega^{k-1}\otimes 1_a \\ 
&&\qquad + n_{a\theta}
\Omega^k\otimes 1_\theta  (\text{if gr (a) =1});\\[3mm]
D(\Omega^k\otimes 1_a) &=& -\disp{ \sum_{b\in \M^{*}_Y(\s)\atop
gr(a)-gr(b)=1}}n_{ab} 
\Omega^k\otimes 1_b;\\[3mm]
D(\Omega^k\otimes 1_\theta) &=&  \disp{\sum_{d\in \M^{*}_Y(\s)\atop
gr(d) = -2}} 
n_{\theta d} \Omega^k\otimes 1_d.
\end{array}
\label{D}
\na
where $n_{ab}, n_{a\theta}$ and $n_{\theta d}$ is the counting
of flowlines  from $a$ to $b$ (if $gr (a) - gr(b) =1$),
 from $a$ to $\theta$ (if $ gr (a) =1$)
and from $\theta$ to $d$ (if $gr(d) =-2$), and $m_{ac}$ (if
$gr(a)-gr(c) =2$) is described as a relative Euler number associated to the
2-dimensional moduli space of flowlines from $a$ to $c$ (Cf. Lemma 5.7 of \cite{MW}).
In the next section, we shall briefly review the construction and
various relations among the coefficients, as 
established in \cite{MW}. These identities ensure that $D^2 =0$.
Notice that, in the complex $CF^{SW}_{*,
U(1)}(Y, \s)$ and in the expression of the differential operator, only
terms with non-negative powers of $\Omega$ are considered. 
We modify the construction as follows.

\begin{Def} Let $CF^{SW, \infty}_{*, U(1)}(Y, \s)$ be the graded
complex generated by
\[
\{\Omega^k\otimes \eta_a, \Omega^k\otimes 1_a, \Omega^k\otimes 1_\theta:
a\in \M^*_Y(\s), k\in \Z.\}
\]
with the grading $gr$ and the differential operator $D$ given
by (\ref{grade}) and (\ref{D}) respectively. Let 
 $CF^{SW, -}_{*, U(1)}(Y, \s)$ be the subcomplex of 
$CF^{SW, \infty}_{*, U(1)}(Y, \s)$, generated by those generators with
negative power of $\Omega$. The quotient complex is
denoted by $CF^{SW, +}_{*, U(1)}(Y, \s)$. Their homologies are denoted
by $HF^{SW, \infty}_{*, U(1)}(Y, \s)$, $HF^{SW, -}_{*, U(1)}(Y, \s)$
and $HF^{SW, +}_{*, U(1)}(Y, \s)$ respectively.  
\end{Def}

The main results in this paper relate these homologies to the
equivariant Seiberg-Witten-Floer homology $HF_{*, U(1)}^{SW}(Y, \s)$
and cohomology $HF_{U(1)}^{SW, *}(Y, \s)$
constructed in \cite{MW} and establish some of their main properties.

\begin{The}  For any rational homology 3-sphere $Y$ with a
$\spinc$ structure $\s \in \spinc (Y)$, these homologies satisfy
the following properties: 
\begin{enumerate}
\item $HF^{SW, \infty}_{*, U(1)} (Y, \s) \cong \Z[\Omega, \Omega^{-1}]$.
\item  $HF^{SW, +}_{*, U(1)} (Y, \s) \cong HF^{SW}_{*, U(1)} (Y, \s)$
where $HF^{SW}_{*, U(1)} (Y, \s)$ is the equivariant
Seiberg-Witten Floer homology for  $(Y, \s)$ constructed in \cite{MW}.
\item  $HF^{SW, -}_{*, U(1)} (Y, \s) \cong HF^{SW, *}_{U(1)} (-Y, \s)$
where $HF^{SW, *}_{U(1)} (-Y, \s)$ is the equivariant
Seiberg-Witten Floer cohomology for  $(-Y, \s)$ constructed in \cite{MW}.
\item There exists a long exact sequence
\ba 
\diagram
\cdots \to HF^{SW, -}_{*, U(1)} (Y, \s) \rto^{l_*} & HF^{SW,
\infty}_{*, U(1)} (Y, \s) 
\rto^{\pi_*} & HF^{SW, +}_{*, U(1)} (Y, \s) \rto^{\delta_*} &
HF^{SW, -}_{*-1, U(1)} (Y, \s) \to \cdots
\enddiagram
\label{exact:main}
\na
relating these homologies. Moreover, $HF^{SW, -}_{*, U(1)} (Y, \s)$,
$HF^{SW, \infty}_{*, U(1)} (Y, \s)$, $HF^{SW, +}_{*, U(1)} (Y, \s)$
and $HF_{red, *}^{SW}(Y, \s) = Coker (\pi_*) \cong Ker (l_{*-1})$ are
all topological invariants of $(Y, \s)$. 
\item There is a $u$-action on $HF^{SW, -}_{*, U(1)} (Y, \s)$,
$HF^{SW, \infty}_{*, U(1)} (Y, \s)$ and $HF^{SW, +}_{*, U(1)} (Y, \s)$
respectively which decreases the degree by two, and is related to
the cutting down moduli spaces of flowlines by a geometric representative
of a degree 2 characteristic form. The long exact sequence
(\ref{exact:main}) is
a long exact sequence of $\Z[u]$-modules.
\item There is a homology group $\widehat{HF}^{SW}_{*} (Y, \s)$, which is
also a topological invariant of $(Y, \s)$,
such that  the following sequence  is exact:
{\small \ba\diagram
\label{exact:hat}
\cdots 
\to \widehat{HF}^{SW}_{*} (Y, \s) \rto & HF^{SW, +}_{*, U(1)} (Y, \s)
  \rto^{u} & HF^{SW, +}_{*-2, U(1)} (Y, \s)\rto &
\widehat{HF}^{SW}_{*-1} (Y, \s) \to\cdots \enddiagram
\na}
and that $\widehat{HF}^{SW} (Y, \s)$ is non-trivial if and only if
$HF^{SW, +}_{*, U(1)} (Y, \s)$ is non-trivial.
\end{enumerate}\label{main-the}
\end{The}

The  $u$-action in the main theorem is induced from 
a $u$-action on the chain complex
\[
u: \qquad CF^{SW, \infty}_{*, U(1)} \to   CF^{SW, \infty}_{*, U(1)},
\]
which decreases the degree by 2. We will show that this $u$-action is
homotopic to the 
obvious $\Omega^{-1}$ action on the chain complex $ CF^{SW,
\infty}_{*, U(1)}$. 
Thus, the induced $u$-action on $HF^{SW, \pm}_{*, U(1)}(Y, \s)$
endows them with $\Z[u]$-module structures.

Let  $\widehat{CF}^{SW}_{*}(Y, \s)$ be the subcomplex of 
$CF^{SW, +}_{*, U(1)} (Y, \s)$ such that 
the following sequence is a short exact sequence of chain complexes:
\[
\diagram
0 \to \widehat{CF}^{SW}_{*}(Y, \s) \rto & CF^{SW, +}_{*, U(1)} (Y, \s)
  \rto^{\Omega^{-1}} & CF^{SW, +}_{*, U(1)} (Y, \s)\to  0
\enddiagram
\]
We can define $\widehat{HF}^{SW}_{*} (Y, \s) $ to be the homology of 
$\widehat{CF}^{SW}_{*}(Y, \s)$.

In recent work \cite{OS1} \cite{OS2}, Ozsv\'ath and Szab\'o introduced
Heegaard Floer invariants $HF^{\pm}_*(Y,\s)$, $HF^{\infty}_*(Y, \s)$,
$\widehat{HF}_*(Y, \s)$, and $HF_{red, *}(Y, \s)$, with exact
sequences relating them. In view of their construction, the result of
Theorem \ref{main-the}, together with the identification of our $HF^{SW,
\infty}_{*, U(1)}(Y, \s)$ and the $HF^{\infty}_*(Y, \s)$ of Ozsv\'ath
and Szab\'o, suggest the following conjecture.

\begin{Conj} For any rational homology 3-sphere $Y$ with a $\spinc$
structure $\s\in \spinc (Y)$, there are isomorphisms
\[\begin{array}{ccc}
HF^{SW, +}_{*, U(1)}(Y, \s) \cong HF^{+}_*(Y, \s), &\qquad &
HF^{SW, -}_{*, U(1)}(Y, \s) \cong HF^{-}_*(Y, \s);\\[3mm]
\widehat{HF}^{SW}_{*}(Y, \s) \cong \widehat{HF}_*(Y, \s), &\qquad &
HF^{SW}_{red, *}(Y, \s) \cong HF_{red, *}(Y, \s).
\end{array}
\]
\end{Conj}

\vskip .2in
\noindent
{\bf Acknowledgments} This research was supported in part  by the Humboldt
Foundation's Sofja Kovalevskaya Award.

\section{Review of equivariant Seiberg-Witten Floer homology}

In this section, we recall some of basic ingredients in the definition
of the   equivariant Seiebrg-Witten Floer homology from \cite{MW} (See
\cite{MW} for all the details).

Let $(Y, \s)$ be a rational homology 3-sphere $Y$ with a $\spinc$
structure $\s\in \spinc (Y)$. For a good pair of metric and perturbation 
(a co-closed imaginary-valued 1-form $\nu$ ) on $Y$, the 3-dimensional 
Seiebrg-Witten equations on $(Y, \s)$
(Cf. \cite{CW} \cite{Fro} \cite{Mar2} \cite{MW}):
\ba\left\{ \begin{array}{l}
*F_A = \sigma (\psi, \psi) + \nu\\
\dirac_A\psi =0,  \\
\end{array} \right.
\label{SW:3d}
\na
for a pair of $\spinc$ connection $A$ and a spinor $\psi$, have only finitely
many irreducible solutions (modulo the gauge transformations), denoted
by $\M^*_Y(\s)$ the set of equivalence classes of irreducible
 solutions to (\ref{SW:3d}),  and  $\theta$ is the 
unique reducible solution (modulo the gauge transformations). Write
$\M_Y(\s) = \M^*_Y(\s) \cup \{\theta\}.$ 

Gauge classes of finite energy solutions to the 4-dimensional Seiebrg-Witten
equations, perturbed as in  \cite{CW} \cite{Fro} \cite{MW},
can be regarded as moduli spaces of flowlines of the Chern-Simons-Dirac
functional on the gauge equivalence classes of $\spinc$ connections
and spinors for $(Y, \s)$. These can be partitioned into 
moduli spaces of flowlines between pairs of critical points from
$\M_Y(\s)$. Each is a smooth oriented  manifold which
can be compactfied to a smooth manifold with corners by adding
broken flowlines that split through intermediate critical points. 

The spectral flow of the Hessian operator of the Chern-Simons-Dirac
functional defines a relative grading
on $\M_Y(\s)$:
\[
gr (\cdot, \cdot): \qquad   \M_Y(\s) \times \M_Y(\s) \rightarrow \Z.\]
In particular, using the unique reducible point $\theta$ in $\M_Y(\s)$, there
is a $\Z$-lifting of the relative grading given by
$gr (a) = gr (a, \theta)$.

Let $a$ be an irreducible monopole in $\M_Y(\s)$, then for any
$b\neq a$ in $ \M_Y(\s)$, the moduli space of flowlines from $a$ to $b$,
denoted by $\M(a, b)$  has dimension 
$gr (a) - gr(b)>0$ (if non-empty). The moduli space of flowlines from
$\theta$ to 
$d\in \M^*_Y(\s)$, denoted by $\M(\theta, d)$ 
has dimension $-gr (d) -1 >0$ (if non-empty). Note that all these
moduli spaces of 
flowlines admit an $\R$-action given by the $\R$-translation of
flowlines: the corresponding quotient spaces are denoted by
$\widehat \M(a, b)$ and $\widehat \M(\theta, d)$, respectively.

For any irreducible critical points $a$ and $c$ in $\M_Y(\s)$
with $gr (a) -gr (c) =2$, we can construct a canonical complex line bundle
over $\M(a, c)$ and a canonical section as follows (see section 5.3
in  \cite{MW}). Choose a base point $(y_0,t_0)$ on $Y\times
\R$, and a complex line $\ell_{y_0}$ in the fiber $W_{y_0}$ of the
spinor bundle $W$ over $y_0\in Y$. We choose $\ell_{y_0}$
so that it does not contain the spinor part $\psi$ of any irreducible
critical point. Since there are only finitely many critical points we
can guarantee such choice exists. Denote the based moduli
space of $\M(a, c)$ by $\M(O_a, O_c)$ as in \cite {MW}, where
$O_a$ and $O_c$ are the $U(1)$-orbits of monopoles 
on the based configuration space. We consider the line bundle
\ba
 {\cal L}_{ac} =\M(O_a,O_c)\times_{U(1)} (W_{y_0}/\ell_{y_0}) \to
\M(a,c) \label{Lac} \na
with a section given by
\ba s([A,\Psi])=([A,\Psi],\Psi(y_0,t_0)). \label{sec:ac} 
\na
For  a generic choice of $(y_0,t_0)$ and $\ell_{y_0}$,
 the section $s$ of (\ref{sec:ac})
has no zeroes on the boundary strata of the compactification of $\M(a, c)$.
This determines a
trivialization of ${\cal L}_{ac}$ away from a compact set in $\M(a,c)$.
The line bundle ${\cal L}_{ac}$ over $\M(a,c)$, with the
trivialization $\varphi$ specified above,
has a well-defined relative Euler class (Cf. Lemma 5.7 in \cite{MW}).

\begin{Def}\label{n:m}
\begin{enumerate}\item For any two irreducible critical 
points $a$ and $b$ in $\M_Y (\s)$ with $gr(a) -gr(b) =1$, we define
$n_{ab}:= \#(\hat\M(a, b))$, the number of flowlines in $\M(a, b)$
counting with orientations. Similarly, for any  $a\in \M_Y (\s)$
with $gr (a) =1$ and any $d\in \M_Y (\s)$ with $gr (d) =-2$, we define
$n_{a\theta} := \#(\hat\M(a, \theta))$ and $n_{\theta d} :=\# (\hat \M
(\theta, b))$, respectively.
\item For any  two irreducible critical 
points $a$ and $c$ in $\M_Y (\s)$ with $gr(a) -gr(c) =1$, we define
$m_{ac}$ to be the relative Euler number of the canonical
line bundle $ {\cal L}_{ac}$  (\ref{Lac}) with the canonical
trivialization given by 
(\ref{sec:ac}).
\end{enumerate}
\end{Def}

The following proposition states various relations satisfied by 
the integers defined in Definition \ref{n:m}, whose proof can be
found in Remark 5.8 of \cite{MW}.

\begin{Pro}
\begin{enumerate} 
\item For any  irreducible critical 
point  $a$  in $\M^*_Y (\s)$ and any critical point $c$ in $\M_Y(\s)$
with $gr (a) -gr(c) =2$,  we have the following identity:
\[
\sum_{b\in \M^*_Y(\s) \atop gr (a) -gr(b) =1} n_{ab}n_{bc} =0.
\]
\item Let $a$ and $d$ be two irreducible critical points with
$gr(a)-gr(d)=3$. Assume that all the critical points $c$ with
$gr(a)> gr(c)> gr (d)$ are irreducible. Then we have the identity
$$ \sum_{c_1: gr (a)-gr (c_1)=1} n_{a,c_1} m_{c_1,d} -
\sum_{c_2:gr (c_2) -gr (d)=1} m_{a,c_2}n_{c_2,d}=0. $$
When $gr(a)=1$ and $gr (d) =-2$, we have the identity
\[
\sum_{c_1:  gr (c_1)=0} n_{a,c_1} m_{c_1,d}
 + n_{a\theta} n_{\theta d}
-\sum_{c_2:gr (c_2) =- 1} m_{a,c_2}n_{c_2,d} =0.
\]
\end{enumerate}
\end{Pro}

With the help of this Proposition, we can check that the
equivariant Seiberg-Witten-Floer complex $CF^{SW}_{*, U(1)}(Y, \s)$
as given in (\ref{CF+}) with the grading and the differential operator
given by (\ref{grade}) and (\ref{D}) is well-defined, and we denote
its homology by $HF^{SW}_{*, U(1)}(Y, \s)$. The
equivariant Seiberg-Witten-Floer cohomology, denoted
by $HF^{SW, *}_{U(1)}(Y, \s)$, is the homology of the dual complex 
$Hom (CF^{SW}_{*, U(1)}(Y, \s), \Z)$. The main result in \cite{MW}
shows that the equivariant Seiberg-Witten Floer homology
$HF^{SW}_{*, U(1)}(Y, \s)$ and 
cohomology $HF^{SW, *}_{U(1)}(Y, \s)$ are topological invariants of
$(Y, \s)$.

\section{Variants of equivariant Seiberg-Witten Floer homology}

As mentioned in the introduction, we will generalize the construction of
the  equivariant Seiberg-Witten Floer homology in several ways. 

First, we denote by $CF^{SW, \infty}_{*, U(1)}(Y, \s)$ the graded complex
generated by
\[
\{\Omega^k\otimes \eta_a, \Omega^k\otimes 1_a, \Omega^k\otimes 1_\theta:
a\in \M^*_Y(\s), k\in \Z.\}
\]

More precisely, for any irreducible critical orbits $O_a$, we set 
\[\begin{array}{lll}
C^\infty_{*, U(1)} (O_a) &=& \Z[\Omega, \Omega^{-1}] \otimes \Omega_0^*(O_a)\\
&:=& \bigoplus_{k\in \Z}\bigl( \Z \Omega^k \otimes \eta_a + 
 \Z\Omega^k \otimes 1_a \bigr)\end{array}
\]
with the grading $gr( \Omega^k \otimes \eta_a) = 2k +gr (a)$ and 
$gr(\Omega^k \otimes 1_a) = 2k+gr(a)+1$, and we set
\[
C^\infty_{*, U(1)} (\theta) = \bigoplus_{k\in \Z} \Z.\Omega^k \otimes 1_\theta
\]
with the grading $gr (\Omega^k \otimes 1_\theta) = 2k$. 

We then consider
\ba
CF^{SW, \infty}_{*, U(1)}(Y, \s) = \bigoplus_{a\in \M_Y(\s)} 
\Z[\Omega, \Omega^{-1}] \otimes \Omega_0^{*-\dim(O_a)}(O_a),
\na
with the grading and the differential operator given by 
(\ref{grade}) and (\ref{D}) respectively. That is, 
$CF^{SW, \infty}_{*, U(1)}(Y, \s)$ is given by
\[
 \begin{array}{lll}
&&\disp{\bigoplus_{a\in \M_Y(\s)}  }C^\infty_{*, U(1)} (O_a)\\
&=& \disp{\bigoplus_{a\in \M^*_Y(\s)}}  C^\infty_{*, U(1)} (O_a) 
\oplus C^\infty_{*, U(1)} (\theta).
\end{array}
\]

\begin{The} \label{HF:infty}
Define $HF^{SW, \infty}_{*, U(1)}(Y, \s)$ to be the
homology of $(CF^{SW, \infty}_{*, U(1)}(Y, \s), D)$. 
Then we have
\[
HF^{SW, \infty}_{*, U(1)}(Y, \s) \cong \Z [\Omega, \Omega^{-1}]. 
 \]
\end{The}
\begin{proof} 
Consider the filtration of $CF^{SW, \infty}_{*, U(1)}(Y, \s)$ according to
the grading of the critical points
$$
 \F_n C^\infty_{*,U(1)}:= \bigoplus_{gr(a) \leq n}
 C_{*,U(1)}^\infty(O_a)
$$
the corresponding spectral sequence $E^r_{kl}$. The filtration is
exhaustive, that is, 
\[
CF^{SW, \infty}_{*, U(1)}(Y, \s)= \bigcup_n \F_n C^\infty_{*,U(1)},
\]
and
\[
\cdots \subset \F_{n-1} C^\infty_{*,U(1)}\subset \F_n C^\infty_{*,U(1)}
\subset \F_{n+1} C^\infty_{*,U(1)} \subset \cdots \subset 
CF^{SW, \infty}_{*, U(1)}(Y, \s).
\] 

Moreover, by the compactness of the moduli space
of critical orbits, the set of indices $gr(a)$ is bounded from above and below,
hence the filtration is bounded. Thus, the spectral sequence
converges to $HF^{SW,\infty}_{*,U(1)}(Y, \s)$.

We compute the $E^0$-term:
\[\begin{array}{lll}
E^0_{kl}&=&\F_k C^\infty_{k+l, U(1)}/ \F_{k-1} C^\infty_{k+l,
U(1)} \\
&= &\disp{\bigoplus_{a\in \M_Y(\s): gr(a)=i \leq k} }
C_{k+l-i,U(1)}^\infty(O_a) /
 \disp{ \bigoplus_{a \in \M_Y(\s): gr( a)=i\leq k-1}
}C_{k+l-i,U(1)}^\infty(O_a) \\  
&= & \disp{\bigoplus_{a \in \M_Y(\s):gr (a)=k}} C_{l,U(1)}^\infty(O_a). 
\end{array}
\]
For $k\neq 0$ this complex is just the direct sum of the separate
complexes $(C_{*,U(1)}^\infty(O_a),\partial_{U(1)})$ on each orbit $O_a$
with $gr(a)=k$:
\ba
\label{acyclic}
 \cdots \to \Z. \Omega\otimes 1_a \stackrel{0}{\to}\Z. \Omega\otimes \eta_a
\stackrel{-1}{\to}\Z. 1\otimes 1_a \stackrel{0}{\to}\Z. 1\otimes \eta_a
\stackrel{-1}{\to} \Z.\Omega^{-1} \otimes 1_a \to \cdots
\na
In the case $k=0$ we have
\[
E^0_{0, l} =  C_{l,U(1)}^\infty(\theta) \oplus \bigoplus_{a\in
\M^*_Y(\s): 
gr(a)=0} C_{l,U(1)}^\infty(O_a),
 \]
which again is a direct sum of the complexes
$(C_{*,U(1)}^\infty(O_a),\partial_{U(1)})$,
here $\partial_{U(1)}$ is the equivariant de Rham differential,
 and of the complex
with generators $\Omega^r\otimes
1_\theta\rangle$ in degree $l=2r$ and trivial differentials.

We then compute the $E^1_{pq}$ term directly: we have
\[ E^1_{kl}= H_{k+l}( E^0_{k,*} ) = \left\{ \begin{array}{ll} \Z. 
\Omega^{r} \otimes 1_\theta  & k =0, l=2r \\[3mm]
0 & k \neq 0, \end{array}\right. \]
since each complex (\ref{acyclic}) is acyclic.
 Thus, the only non-trivial $E^1$-terms are of the form $E^1_{0l}=\Z. 
\Omega^{r} \otimes 1_\theta$, $l=2r$, with trivial differentials,
so that the spectral sequence collapses and we obtain the result.
\end{proof}

\subsection{Long exact sequence}

\begin{Def} Let $CF^{SW, -}_{*, U(1)}(Y, \s)$ be the subcomplex of
$CF^{SW, \infty}_{*, U(1)}(Y, \s)$, generated by 
\[
\{\Omega^k\otimes \eta_a, \Omega^k\otimes 1_a , \Omega^k\otimes 1_\theta:
a\in \M^*_Y(\s), k\in \Z \ and \ k < 0\},
\]
whose homology groups are denoted by $HF^{SW, -}_{*, U(1)}(Y, \s)$.
The quotient complex is denoted by $CF^{SW, +}_{*, U(1)}(Y, \s)$,
with the homology groups denoted by $HF^{SW, +}_{*, U(1)}(Y, \s)$.
\end{Def}

\begin{The}\label{exact:sequence}
\begin{enumerate} 
\item $HF^{SW, +}_{*, U(1)} (Y, \s) \cong HF^{SW}_{*, U(1)}(Y, \s)$,
where $HF^{SW}_{*, U(1)}(Y, \s)$ is the equivariant
Seiberg-Witten-Floer homology 
defined in \cite{MW}.
\item There is an exact sequence of $\Z$-modules which relates these variants
of equivariant Seiberg-Witten-Floer homologies:
\[
\diagram
\cdots \to HF^{SW, -}_{*, U(1)} (Y, \s) \rto^{l_*} & HF^{SW,
\infty}_{*, U(1)} (Y, \s) 
\rto^{\pi_*} & HF^{SW, +}_{*, U(1)} (Y, \s) \rto^{\delta_*} &
HF^{SW, -}_{*-1, U(1)} (Y, \s) \to \cdots
\enddiagram
\]
\end{enumerate}
\end{The}
\begin{proof} It is easy to see that $CF^{SW, +}_{*, U(1)}(Y, \s) 
= CF^{SW}_{*, U(1)}(Y, \s)$, with the same grading and differentials,
hence $HF^{SW, +}_{*, U(1)} (Y, \s) \cong HF^{SW}_{*, U(1)}(Y, \s)$.
The long exact sequence in homology is induced by the short exact
sequence of chain complexes: 
\[
0 \to CF^{SW, -}_{*, U(1)}(Y, \s) \to  CF^{SW,\infty}_{*, U(1)}(Y, \s)
\to  CF^{SW,+}_{*, U(1)}(Y, \s)\to 0.
\]
\end{proof}

From the above long exact sequence, we can define 
\ba
\begin{array}{lll}
HF^{SW}_{red, *}(Y, \s) & =& 
Coker (\pi_*) \cong HF^{SW, +}_{*, U(1)}(Y, \s)/Ker (\delta_*) \\
&\cong & Im (\delta_*) \cong Ker (l_{*-1}).
\end{array}
\label{HF_red}
\na

\subsection{The spectral sequence for $HF_{*,U(1)}^{SW,+}(Y,\s)$} 
 
We consider again the filtration by index of critical orbits, 
$$ \cF_n C_{*,U(1)}^+ := \bigoplus_{\gr(a)\leq n} 
C^+_{*,U(1)}(O_a), $$ for 
$$ C^+_{*,U(1)}(O_a) = \Z [\Omega] \otimes \Omega_0^{*-\dim(O_a)}(O_a). $$ 
 
We have 
$$ 
\begin{array}{rl} E^0_{kl} = & \cF_k C^+_{k+l, U(1)} /\cF_{k-1}
C^+_{k+l, U(1)}  
\\[3mm] = & \bigoplus_{\gr(a)=k} C^+_{l,U(1)}(O_a). \end{array} 
$$ 
This is a direct sum of the complexes 
\begin{equation} 
\label{orbit-U(1)} 
\cdots \stackrel{-1}{\to} \Z.\Omega \otimes 1_a \stackrel{0}{\to} 
\Z.\Omega \otimes \eta_a \stackrel{-1}{\to} \Z.1 \otimes 1_a 
\stackrel{0}{\to}\Z. 1 \otimes \eta_a \to 0, 
\end{equation} 
over each orbit $O_a\cong S^1$ and, in the case $k=0$, the complex 
with generators $\Omega^r \otimes 1_\theta$ in degree $l=2r$, and 
trivial differentials. 
 
Thus, we obtain that $E^1_{pq} = H_{p+q} (E^0_{p*})$ is of the 
form 
$$ 
E^1_{pq} = \left\{ \begin{array}{ll} 0 & q>0 \\[2mm] 
\Z . 1\otimes \eta_a  & q=0, \gr(a)=p \end{array} \right. 
$$ 
for $p\neq 0$, and 
$$ 
E^1_{0q} = \left\{ \begin{array}{ll} \Z. \Omega^r \otimes 1_\theta  &
q=2r >0 \\[2mm]  
\Z . 1\otimes \eta_a  \oplus \Z. 1\otimes 1_\theta & q=0, 
\gr(a)=0. 
\end{array} \right. 
$$ 
 
The differential $d^1: E^1_{p,q}\to E^1_{p-1,q}$ is of the form 
$$ \begin{array}{rl} d^1(1\otimes \eta_a) = & n_{ab} 1\otimes 
\eta_b \\[2mm] 
& + n_{a\theta} 1\otimes 1_\theta \,\, \text{ (if $\gr(a)=1$) } 
\end{array} $$ 
 
Thus, we obtain 
$$ 
E^2_{pq} = \left\{ 
\begin{array}{ll} HF_p^{SW}(Y,\s) & p\neq 0, 
q=0 \\[3mm] 
Ker(\Delta_1) & p=1, q=0 \\[3mm] 
HF_0^{SW}(Y,\s) \oplus T_0 & p=0, q=0 \\[3mm] 
\Z.\Omega^r\otimes 1_\theta & p=0, q=2r >0. 
\end{array} 
\right. 
$$ 
Here $HF_*^{SW}(Y,\s)$ denotes the non-equivariant (metric and 
perturbation dependent) Seiberg--Witten Floer homology. This is 
the homology of the complex with generators $1\otimes \eta_a$ in 
degree $\gr(a)$ and boundary coefficients $n_{ab}$ for 
$\gr(a)-\gr(b)=1$. We also denoted by $\Delta_1$ the map 
$$ \Delta_1: HF_1^{SW}(Y,\s) \to \Z. 1\otimes 1_\theta, $$ 
$$ 
\Delta_1 (\sum x_a 1\otimes \eta_a) = \sum x_a n_{a\theta} 
1\otimes 1_\theta, 
$$ 
where the coefficients $x_a$ satisfy $\sum x_a n_{ab}=0$. Finally, 
the term $T_0$ denotes the term 
$$ T_0 = \Z. 1\otimes 1_\theta / Im(\Delta_1). $$ 
 
Notice then that the boundary $d^2: E^2_{p,q} \to E^2_{p-2,q+1}$ 
is trivial, hence the $E^3_{p,q}$ terms are disposed as in the 
diagram: {\small 
$$ 
\spreaddiagramrows{-1pc} \spreaddiagramcolumns{-1.75pc} 
\diagram 
\cdots & 0 & 0 & 0 & 0 & \Z . \Omega^2\otimes 1_\theta & 0 & 0 &\cdots 
\\ 
\cdots & 0 & 0 & 0 & 0 & 0 & 0  & 0& \cdots 
\\ 
\cdots & 0 & 0 & 0 & 0 & \Z . \Omega \otimes 1_\theta & 0 & 0 & \cdots 
\\ 
\cdots & 0 & 0 & 0 & 0 & 0 & 0 & 0 & \cdots 
\\ 
\cdots & HF_{4}^{SW} & HF_3^{SW} \xto[-2,3]^{d^3} & HF_2^{SW} & 
Ker(\Delta_1) & HF_0^{SW} \oplus T_0 & HF_{-1}^{SW} & HF_{-2}^{SW} & \cdots 
\enddiagram 
$$ } 
 
The differential $d^3: E^3_{p,q} \to E^3_{p-3,q+2}$ is given by 
the expression 
\begin{equation}\label{d3} 
 d^3 ( [\sum x_a 1\otimes \eta_a] ) = \sum x_a m_{ac} n_{c\theta} 
\Omega \otimes 1_\theta, 
\end{equation} 
for $\gr(a)-\gr(c)=2$. The expression is obtained by considering 
the unique choice of a representative of the class $[\sum x_a 
1\otimes \eta_a]$ in $E^3_{p,q}$ whose boundary (\ref{D}) 
defines a class in $E^3_{p-3,q+2}$. 
 
The differential $d^4: E^4_{p,q} \to E^4_{p-4,q+3}$ is again 
trivial, and we obtain the $E^5_{pq}$ of the form {\small 
$$ 
\spreaddiagramrows{-1pc} \spreaddiagramcolumns{-1.75pc} \diagram 
\cdots & 0 & 0 & 0 & 0 & \Z . \Omega^2\otimes 1_\theta & 0 & 0  &\cdots 
\\ 
\cdots & 0 & 0 & 0 & 0 & 0 & 0  & 0  &\cdots 
\\ 
\cdots & 0 & 0 & 0 & 0 & T_1 & 0 &0  & \cdots 
\\ 
\cdots & 0 & 0 & 0 & 0 & 0 & 0  & 0  &\cdots 
\\ 
HF_{5}^{SW} \xto[-4,5]^{d^5} & HF_{4}^{SW} & Ker(\Delta_3)  & 
HF_2^{SW} & Ker(\Delta_1) & HF_0^{SW} \oplus T_0 & HF_{-1}^{SW} &  HF_{-2}^{SW} &
\cdots 
\enddiagram 
$$ } 
where again we denote by $T_1$ the  term 
$$ T_1 = \Z . \Omega \otimes 1_\theta / Im (\Delta_3). $$ 
 
Thus, by iterating the process, we observe that all the 
differentials $d^{2k}: E^{2k}_{p,q} \to E^{2k}_{p-2k,q+2k+1}$ are 
trivial and the differentials $d^{2k+1}: E^{2k+1}_{p,q} \to 
E^{2k+1}_{p-2k-1, q+2k}$ consists of one map for $p=2k+1$, $q=0$: 
$$ \Delta_{2k+1} : HF_{2k+1}^{SW} \to \Z . \Omega^k \otimes 1_\theta, $$ 
induced by 
$$ 
\Delta_{2k+1} ( \sum x_a 1\otimes \eta_a) = \sum x_a m_{a 
a_{2k-1}} m_{a_{2k-1} a_{2k-3} } \cdots m_{a_3 a_1 } n_{a_1\theta 
} \Omega^k \otimes 1_\theta. 
$$ 
Here we have $\gr(a)=2k+1$ and $\gr (a_r)=r$. Notice that these 
maps agree with the morphism $\Delta_*$, which is obtained in 
\cite{MW} as the connecting homomorphism in the long exact 
sequence relating equivariant and non-equivariant Seiberg--Witten 
Floer homologies. 
 
We thus obtain the following structure theorem for equivariant 
Seiberg--Witten Floer homology. 
 
\begin{The} 
The equivariant Seiberg--Witten Floer homology $HF_{*,U(1)}^{SW,+}(Y,\s)$
has the form 
$$ HF_{*,U(1)}^{SW,+}(Y,\s)= \left\{ \begin{array}{ll}
Ker(\Delta_{2k+1}) & *=2k+1 >0 \\[3mm]  
HF_{2k}^{SW} (Y,\s) \oplus T_{k} & *=2k \geq 0 \\[3mm] 
HF_*^{SW}(Y,\s) & * < 0 \end{array} \right. 
$$ 
where $T_{k}$ is the term 
$$ T_{k} = \Z . \Omega^k \otimes 1_\theta / Im (\Delta_{2k+1}). 
$$ 
\end{The}

This result refines the long exact sequence obtained in \cite{MW}: 
$$
\spreaddiagramrows{-1pc}
\spreaddiagramcolumns{-1pc}
\diagram
HF^{SW}_{*,U(1)}(Y,\s) \rto^{i_*} & HF^{SW}_* (Y,\s,g,\nu) 
\dlto^{\Delta_*} \\
\Z [\Omega ]\uto^{j_*} 
\enddiagram
$$

Similar results can be obtained for $HF_{*,U(1)}^{SW,-}(Y,\s)$.

\subsection{Topological invariance}

Note that the definitions of these homologies depend on the
Seiberg-Witten equations, which use the metric and perturbation
on $(Y, \s)$.
By the result of \cite{MW}, we know that
$HF^{SW, +}_{*, U(1)} (Y, \s) \cong HF^{SW}_{*, U(1)}(Y, \s)$ is 
a topological invariant of $(Y, \s)$, we first recall this topological
invariance as stated in Theorem 6.1 \cite{MW}.

\begin{The}\label{topo:inv}
 (Theorem 6.1 \cite{MW})  Let $(Y, \s)$ be a
rational homology sphere with a $\spinc$ structure.
Suppose given two metrics $g_0$ and $g_1$ on $Y$ and perturbations
$\nu_0$ and $\nu_1$ such that
$Ker(\dirac^{g_0}_{\nu_0})=Ker(\dirac^{g_1}_{\nu_1})=0$, so that  
the corresponding monopole moduli spaces
 $\M_Y(\s, g_0, \nu_0)$ and  $\M_Y(\s, g_1, \nu_1)$ consist of
finitely many isolated points. 
Then there exists   an isomorphism between 
the equivariant Seiberg-Witten Floer homologies $HF^{SW}_{*, U(1)}(Y,
\s, g_0, \nu_0)$ 
and  $HF^{SW}_{*, U(1)}(Y, \s, g_1, \nu_1)$, with a degree shift given
by the spectral flow of the
Dirac operator $\dirac^{g_t}_{\nu_t}$ along a path of metrics and
perturbations 
connecting $(g_0, \nu_0)$ and $(g_1, \nu_1)$. That is, 
if  the complex spectral flow along the path
$(g_t,\nu_t)$ is denoted by $SF_\C(\dirac^{g_t}_{\nu_t})$, then for any $k\in \Z$, 
\[
HF^{SW}_{k, U(1)}(Y, \s, g_0,\nu_0) \cong HF^{SW}
_{k+2SF_\C(\dirac^{g_t}_{\nu_t}), 
U(1)}(Y, \s, g_1,\nu_1).
\]
\end{The}

   From  Theorem \ref{HF:infty}, we know that 
\[
HF^{SW, \infty}_{*, U(1)} (Y, \s) \cong \Z[\Omega, \Omega^{-1}]
\]
is independent of $(Y, \s)$, up to a degree shift as given in Theorem
\ref{topo:inv}. Thus, applying the five lemma to the 
long exact sequence in Theorem \ref{exact:sequence}, we obtain  
that  $HF^{SW, -}_{*, U(1)} (Y, \s)$ and $HF^{SW}_{red, *}(Y, \s) $ are also
topological invariants of $(Y, \s)$.

\begin{The} $HF^{SW, -}_{*, U(1)} (Y, \s)$ and $HF^{SW}_{red, *}(Y,
\s) $ are  
topological invariants of $(Y, \s)$, in the sense that,
given any two metrics $g_0$ and $g_1$ on $Y$ and perturbations
$\nu_0$ and $\nu_1$, with  
$Ker(\dirac^{g_0}_{\nu_0})=Ker(\dirac^{g_1}_{\nu_1})=0$, 
there exist isomorphisms 
\[\begin{array}{c}
HF^{SW, -}_{k, U(1)}(Y, \s, g_0,\nu_0)
\cong 
 HF^{SW, -}_{k+2SF_\C(\dirac^{g_t}_{\nu_t}),
U(1)}(Y, \s, g_1,\nu_1)\\
HF^{SW}_{red, k} (Y, \s,  g_0,\nu_0)  \cong HF^{SW}_{red,
k+ 2SF_\C(\dirac^{g_t}_{\nu_t})} (Y, \s,  g_1,\nu_1). \\
\end{array}
\]
Here $SF_\C(\dirac^{g_t}_{\nu_t})$ denotes the complex spectral flow
of the Dirac operator $\dirac^{g_t}_{\nu_t}$ along the path
$(g_t,\nu_t)$.
\end{The}

\section{Properties of equivariant Seiberg-Witten Floer homologies}

In this section, we briefly discuss some of the algebraic structures
and properties 
of the equivariant Seiberg-Witten Floer homologies defined in 
the previous section. 

Note that for any irreducible critical points $a$ and $b$ in
$\M^*_Y(\s)$, the associated integer $m_{ac}$ is the counting of
points in the geometric representative of the relative first
Chern class of the canonical line bundle 
(\ref{Lac}) over $\M(a, c)$, we can apply this
fact to define a $u$-action on the chain complex $CF^{SW, \infty}_{*,
U(1)}(Y, \s)$ 
\[
u:\qquad CF^{SW, \infty}_{*, U(1)}(Y, \s)
 \longrightarrow CF^{SW, \infty}_{*, U(1)}(Y, \s)
\]
which decreases the grading by two.  The action is given in terms of
its actions on generators as follows:
\ba
\begin{array}{lll}
u(\Omega^n \otimes \eta_a) &= &\disp{\sum_{c\in\M^*(Y, \s)\atop gr(a)-gr(c) =2}}
m_{ac} \Omega^n \otimes \eta_c.\\
u(\Omega^n \otimes  1_a) &=& \left\{
  \begin{array}{ll}
\disp{\sum_{c\in\M^*(Y, \s)\atop gr(a)-gr(c) =2}}
m_{ac} \Omega^n \otimes 1_c \qquad & \text{if $gr(a) \neq 1$}\\
\disp{\sum_{c\in\M^*(Y, \s)\atop  gr(c) =-1}}
m_{ac} \Omega^n \otimes 1_c + n_{a\theta}\Omega^n\otimes 1_\theta
\qquad & \text{if $gr(a) =1$}\end{array}\right.\\
u(\Omega^n\otimes 1_\theta)& =& \disp{\sum_{d\in \M^*_Y(\s)\atop gr(d) =-2}}
n_{\theta d} \Omega^n\times \eta_d + \Omega^{n-1}\otimes 1_\theta.
\end{array}
\label{u-action}
\na

\begin{Pro} \label{homotopy}
 The u-action defined (\ref{u-action}) on the chain complex
$ CF^{SW, \infty}_{*, U(1)}(Y, \s)$ is homotopic to
the $\Omega^{-1}$-action acting on $ CF^{SW, \infty}_{*, U(1)}(Y, \s)$. 
The induced actions on 
$CF^{SW, \pm}_{*, U(1)}(Y, \s)$ define $\Z[u]$-module structures
on $HF^{SW,^\pm}_{*, U(1)}(Y, \s)$.
\end{Pro}
\begin{proof}
Define $H:  CF^{SW, \infty}_{*, U(1)}(Y, \s)
 \longrightarrow CF^{SW, \infty}_{*, U(1)}(Y, \s)$ by its actions
on the generators as follows:
\[\begin{array}{ll}
&H(\Omega^n \otimes \eta_a) =0,\\
& H(\Omega^n \otimes  1_a) = \Omega^n \otimes \eta_a,\\
& H( \Omega^n\otimes 1_\theta) = 0.\end{array}
\]
Then it is a direct calculation to show that we have:
$$ (u-\Omega^{-1}) (\Omega^k\otimes \eta_a)= m_{ac} \Omega^k \otimes
\eta_c - \Omega^{k-1} \otimes \eta_a = (DH+HD)(\Omega^k\otimes \eta_a) $$
$$ (u-\Omega^{-1}) (\Omega^k\otimes 1_a)= m_{ac} \Omega^k \otimes 1_c
- \Omega^{n-1} \otimes 1_a  \left( + n_{a\theta} \Omega^n \otimes
1_\theta \, \text{ if } \gr(a)=1 \right)= (DH+HD)(\Omega^k\otimes 1_a), $$
$$ (u-\Omega^{-1}) (\Omega^k\otimes 1_\theta) = n_{\theta d} \Omega^n
\otimes \eta_d = (DH+HD) (\Omega^k\otimes 1_\theta). $$ 

Thus the claim follows using the chain homotopy $u-\Omega^{-1} = D\circ
H + H\circ D$. 

\end{proof}

Thus, on the homological level, we can identify the $u$-action with
the induced 
$\Omega^{-1}$ action on various homologies. In particular, we see
that there is a subcomplex $\widehat{CF}^{SW}_*(Y, \s)$ of 
$CF^{SW,+}_{*, U(1)}(Y, \s)$ such that the following short
exact sequence of chain complexes holds:
\ba\diagram
0\to \widehat{CF}^{SW}_*(Y, \s) \rto & 
CF^{SW,+}_{*, U(1)}(Y, \s) \rto^{\Omega^{-1}} & CF^{SW,+}_{*, U(1)}(Y, \s)
\to 0.
\enddiagram
\label{exact:hatCF}
\na

\begin{Pro}
Let $\widehat{HF}^{SW}_*(Y, \s)$ be the homology of
$\widehat{CF}^{SW}_*(Y, \s)$, then 
$\widehat{HF}^{SW}_*(Y, \s)$ is also a topological invariant of $(Y,
\s)$, and it is 
determined by the following long exact sequence
\[\diagram
\cdots   
\to \widehat{HF}^{SW}_{*} (Y, \s) \rto & HF^{SW, +}_{*, U(1)} (Y, \s)
  \rto^{u} & HF^{SW, +}_{*-2, U(1)} (Y, \s)\rto &
\widehat{HF}^{SW}_{*-1} (Y, \s) \to\cdots. \enddiagram
\]
Moreover, $\widehat{HF}^{SW} (Y, \s)$ is non-trivial if and only if
$HF^{SW, +}_{*, U(1)} (Y, \s)$ is non-trivial.
\end{Pro}
\begin{proof}
The long exact sequence follows from the short
exact sequence of chain complexes (\ref{exact:hatCF}) and Proposition
\ref{homotopy}. This long exact sequence  
implies that $\widehat{HF}^{SW}_*(Y, \s)$ is also a topological
invariant of $(Y, \s)$. 

Note that, from the compactness of $\M_Y(\s)$, we see that
each element in $HF^{SW, +}_{*, U(1)} (Y, \s)$ can be annihilated by
a sufficiently 
large power of  $\Omega^{-1}$. Hence, $u$ is an isomorphism on
 $HF^{SW, +}_{*, U(1)} (Y, \s)$ if and only if 
$HF^{SW, +}_{*, U(1)} (Y, \s)$ is trivial. Then the last claim follows from
this observation and the long exact sequence.
\end{proof}

If we think of the set of $\spinc$ structures on $Y$ as the set
of equivalence classes of nowhere vanishing vector fields on $Y$ (Cf.\cite{Tur}),
then there is a natural bijection between $\spinc (Y)$ and $\spinc (-Y)$ where
$-Y$ is the same $Y$ with the opposite orientation.

\begin{The}\label{duality}
Let $(Y, \s)$ be a rational homology 3-sphere with a $\spinc$ structure $\s$,
and $(-Y, \s)$ denote $Y$ with the opposite orientation and the
corresponding $\spinc$ structure. Then there is a natural isomorphism
\[
HF^{SW, *}_{U(1)}(Y, \s) \cong HF^{SW, -}_{*, U(1)}(-Y, \s) 
\]
where $HF^{SW, *}_{U(1)}(Y, \s)$ is the equivariant Seiberg-Witten-Floer
cohomology defined in \cite{MW}.
\end{The}
\begin{proof}
Note that $HF^{SW, *}_{U(1)}(Y, \s)$ is the homology of the dual complex
$Hom(CF^{SW, +}_{*, U(1)}(Y, \s), \Z)$. We start to construct a natural 
pairing 
\ba
\label{pairing}
 \la \cdot, \cdot \ra: \qquad CF^{SW, \infty}_{*, U(1)}(Y, \s) \times 
CF^{SW, \infty}_{*, U(1)}(-Y, \s) \longrightarrow \Z
\na
which satisfies 
\ba
\la D_Y(\xi_1), \xi_2\ra =\la \xi_1, D_{-Y}(\xi_2)\ra,
\qquad \la \Omega^{-1}(\xi_1), \xi_2\ra =\la \xi_1, \Omega^{-1}(\xi_2)\ra.
\label{pair:D}
\na
for any element $\xi_1 \in CF^{SW, \infty}_{*, U(1)}(Y, \s)$ and
any element $\xi_2 \in CF^{SW, \infty}_{*, U(1)}(-Y, \s)$.

Then we will show that the above pairing is non-degenerate when
restricted to
$CF^{SW, +}_{*, U(1)}(Y, \s) \times 
CF^{SW, -}_{*, U(1)}(-Y, \s)$.

From the nature of the Seiberg-Witten equations, we see that
there is an identification 
\[
\M_Y(\s) \to \M_{-Y}(\s)
\]
for a good pair of metric and perturbation on $(Y, \s)$ and
the corresponding metric and perturbation on $(-Y, \s)$.
Then  the relative gradings with respect to the unique reducible
monopole in $\M_Y(\s)$ and $ \M_{-Y}(\s)$ respectively, satisfies
\[
gr_{-Y} (a^-) = - gr_Y (a) -1,
\]
where $a^-$ is the element in $\M^*_{-Y}(\s)$ corresponding to $a\in 
\M^*_Y(\s) $, we assume that $gr_Y (\theta) = gr_{-Y}(\theta^-)$.
Moreover, there is an natural identification between 
the moduli spaces of flowlines for $(Y, \s)$ and $(-Y, \s)$, that is,
\[
\M_{Y\times \R}(a, b) \cong \M_{-Y\times \R}(b^-, a^-).\]

Now we define the pairing on $ CF^{SW, \infty}_{*, U(1)}(Y, \s) \times 
CF^{SW, \infty}_{*, U(1)}(-Y, \s) $
such that the following pairings are the only non-trivial pairings:
\[
\la \Omega^n\otimes \eta_a, \Omega^{-n-1} \otimes 1_{a^-}\ra =1\]
\[
\la \Omega^n\otimes 1_a, \Omega^{-n-1} \otimes \eta_{a^-}\ra =1\]
\[ \la \Omega^n\otimes 1_\theta, \Omega^{-n-1} \otimes 1_{\theta^-}\ra =1.
\]

It is a direct calculation to show that this pairing
satisfies the relation (\ref{pair:D}) and the restriction
of this pairing to $CF^{SW, +}_{*, U(1)}(Y, \s) \times 
CF^{SW, -}_{*, U(1)}(-Y, \s) $ is non-degenerate. Then the claim follows from
the definition.
\end{proof}

Let $\widehat{HF}^{SW, *}(Y, \s)$
and $HF^{SW, *}_{\pm, U(1)}(Y, \s)$ denote the
homology groups of the dual complexes 
$Hom(\widehat{CF}^{SW}_*(Y, \s), \Z)$ and
 $Hom(CF^{SW, \pm}_{*, U(1)}(Y, \s), \Z)$ of 
$\widehat{CF}^{SW}_*(Y, \s)$
and $CF^{SW, \pm}_{*, U(1)}(Y, \s)$
respectively. From the proof the above Theorem \ref{duality}, we actually establish
the following duality between these homologies.

\begin{The} For any rational homology 3-sphere $Y$ with a $spinc$ structure
$\s$, there exist natural isomorphisms
\ba
\label{duality:2}
\widehat{HF}^{SW, *}(Y, \s) \cong \widehat{HF}^{SW}_*(-Y, \s),
\qquad HF^{SW, *}_{\pm, U(1)}(Y, \s) \cong 
 HF^{SW, \mp }_{*, U(1)}(-Y, \s).
\na
\end{The}

\vspace{1cm}

\vskip .3in
 
\noindent {\bf Matilde Marcolli} and {\bf Bai-Ling Wang},
\par\noindent Max--Planck--Institut f\"ur Mathematik,
\par\noindent Vivatsgasse 7, D-53111 Bonn,
Germany.
 
\smallskip
\noindent marcolli\@@mpim-bonn.mpg.de
\par \noindent
bwang\@@mpim-bonn.mpg.de


\begin{thebibliography}{99}
 
 

\bibitem{AB} D.M. Austin, P.J. Braam, {\em Equivariant Floer theory
    and gluing Donaldson polynomials}, Topology 35 (1996), No.1,
  167-200.
 
\bibitem{CW} A.L. Carey,  B.L. Wang, {\em
          Seiberg-Witten-Floer homology and Gluing formulae}, preprint.


\bibitem{Fro} K.A. Fr\/oyshov, {\em The Seiberg-Witten equations and
    four-manifolds with boundary}, Math. Res. Lett. 3 (1996), N.3,
  373--390.

\bibitem{Mar2} M. Marcolli, {\em Seiberg--Witten gauge theory}, Texts
and Readings in Mathematics, 17. Hindustan Book Agency, New Delhi,
1999.

\bibitem{MW} M.Marcolli and B.L.Wang, {\em Equivariant
Seiberg--Witten Floer homology}, Commun. Anal. Geom. Vol.9 N.3
(2001) 451-639.

 
\bibitem{MW:SW=CW}  M. Marcolli, B.L. Wang, {\em
   Seiberg-Witten invariant  and Casson-Walker invariant for any rational
   homology 3-sphere}, Geometriae Dedicata.  Vol.91 (1): 45-58, April 2002.

\bibitem{OS1} P. Ozsvath, Z. Szabo, {\em 
  Holomorphic disks and topological invariants for closed three-manifolds},
  math.SG/0101206.

\bibitem{OS2} P. Ozsvath, Z. Szabo, {\em Holomorphic disks and three-manifold invariants: properties and applications}, math.SG/0105202.

\bibitem{Tur} V. Turaev, {\em Torsion invariants of $\spinc$ structures
on 3-manifolds}, Math. Research Letters, 4, 679-695, 1997.

\end{thebibliography}
\end{document}